\documentclass[times]{my}

\begin{document}

% Enter full title and short title for running headers
\title{Approximation of the Birth Process}
\shorttitle{Approximation of the Birth Process}

% Author name(s)
\author{Mariia Nosova\affil{1}}
% Abbreviated author name for running headers
\abbrevauthor{M. Nosova}
% Abbreviated author name for first page header
\headabbrevauthor{Nosova, M}

\address{%
\affilnum{1}Tomsk University of Control Systems and Radioelectronics, Tomsk, Russia}

% Address / e-mail address of corresponding author
\correspdetails{nosovamgm@gmail.com}

% Received/revised/accepted dates will be entered by the publisher during production of an accepted paper. Please do not edit these placeholders for submission.
\received{}
\revised{}
\accepted{}

% Enter details of editor communicating this article
\communicated{ }

\begin{abstract}

In this article, queuing systems with an unlimited number of devices with an incoming nonstationary Poisson flow and a random flow controlled by a Markov chain are investigated. The inexpediency of approximation of the birth process by Poisson flows in long-term forecasting is shown.
\end{abstract}

\maketitle

\section{Introduction}

In demography, the birth process is often approximated by the Poisson flow \cite{bibid1}, for example in \cite{bibid2, bibid3, bibid4, bibid5, bibid6, bibid7}. Such a model seems to be quite adequate to the actual situation with short-term forecasting for periods of not more than 10-15 years. However, for medium-term forecasts of the order of 20-50 years, it is necessary to consider more adequate mathematical models of the fertility process taking into account the time variation in the number of women of reproductive age, and taking into account the distribution of the age group sizes of those women since fertility rates vary substantially within the reproductive age (15-49 years old) depending on the age of women. 

As such, the mathematical models of the birth process, random flows with random intensity (Cox fluxes) or double stochastic flows, or in other words random streams controlled by random processes among which the most popular MMR flows or MAP streams, can apply \cite{bibid9}. 

In this paper, we consider the simplest mathematical model of a random flow controlled by a Markov chain. Such a model is defined as the flow of customers in an autonomous queuing system with an unlimited number of devices and an average value \textit {m}(\textit {t}) of the number of occupied devices.

And we also consider a queuing system with an unlimited number of devices, the input of which is a non-stationary Poisson flow with intensity $\lambda$(\textit{t}) chosen from the condition of coincidence of the average values of the number of devices occupied in the autonomous system and in a system with a Poisson incoming flow.

We show that the characteristics of such queuing systems differ significantly, which suggests that the Poisson flow approximation is inappropriate for long-term forecasting \cite{bibid7, bibid8}.  

\section{Mathematical model of a Markov autonomous queuing system} 
Let us consider a Markov autonomous queuing system with an unlimited number of devices (Figure 1).

\begin{figure}[h]
\centering
\includegraphics[width=8.493cm,height=6.165cm]{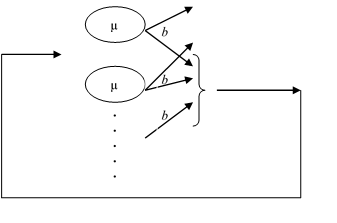}
\caption{Markov autonomous queuing system}
\label{figure: header text example}
\end{figure}

An incoming customer takes any free device. The maintenance time has an exponential distribution with the parameter  $\mu$. During the maintenance time the customer with intensity \textit{b} generates new orders that occupy other free devices that is with the probability \textit{b}$\Delta$\textit{t}+\textit{o}($\Delta$\textit{t}) for an infinitesimal time interval (\textit{t}, \textit{t}+$\Delta$\textit{t}) the service that is served generates a new one. After completing maintenance on the device, the customer leaves the system. The times for servicing various customers and the procedure for generating new applications are stochastically independent. It is obvious that the number of devices \textit{i}(\textit{t}) occupied in the system under consideration is a Markov chain therefore the flow of applications arriving at the devices is a random flow controlled by the Markov chain.

\section{Research of a Markov autonomous queuing system} 

Since the process \textit{i}(\textit{t})\textit{ }is a Markov chain then for its probability distribution 
\[P(i,t)=P\{i(t)=i\}\] 
we write equation
\[P\left(i,t+\Delta t\right)=P\left(i,t\right)\left(1-i\Delta t\right)\left(1-ib\Delta t \right)+P\left(i-1,t\right)\left(i-1\right)b\Delta t +P\left(i+1,t\right)\left(i+1\right)\Delta t+o(\Delta t),\]
and obtain a system of Kolmogorov differential equations
\begin{equation} \label{GrindEQ__1_} 
\frac{\partial P(i,t)}{\partial t}=-i(\mu+b)P(i,t)+(i-1)bP(i-1,t)+(i+1)\mu P(i+1,t).                     
\end{equation} 
Denoting the characteristic function of the number of occupied devices as
\[H(u,t)=Me^{jui(t)}=\sum^\infty_{i=0}{e^{jui}}P(i,t),              \] 
from \eqref{GrindEQ__1_} we obtain the equation for \textit{H}(\textit{u},\textit{t}) in the form
\begin{equation} \label{GrindEQ__2_} 
\frac{\partial H(u,t)}{\partial t}+j\left\{(e^{ju}-1)b+(e^{-ju}-1) \mu \right\}\frac{\partial H(u,t)}{\partial u}=0.                                        
\end{equation} 
We find the solution \textit{H}(\textit{u},\textit{t}) by the method of characteristics. Equation for characteristics
\begin{equation} \label{GrindEQ__3_} 
\frac{dt}{1}=\frac{du}{j\left\{(e^{ju}-1)b+(e^{-ju}-1)\mu \right\}} 
\end{equation} 
is an equation with separated variables whose solution is determined by integrating the left and right parts of it. The first integral of equation \eqref{GrindEQ__3_} can be written in the form
\[C=e^{bt}{\left(\frac{e^{ju}-1}{e^{ju}-\frac{\mu}{b}}\right)}^{\frac{b}{b-\mu}},\] 
therefore, the general solution of equation \eqref{GrindEQ__2_} can be written as follows
\begin{equation} \label{GrindEQ__4_} 
H(u,t)=\phi \left(e^{bt}{\left(\frac{e^{ju}-1}{e^{ju}-\frac{\mu}{b}}\right)}^{\frac{b}{b-\mu}}\right),                                                         
\end{equation} 
where $\phi (.)$ is an arbitrary differentiable function.

We assume that at the initial instant of time \textit{t=}0 in the considered queuing system \textit{N} instruments is occupied that is for the equation \eqref{GrindEQ__2_} an initial condition $H\left(u,0\right)=e^{juN}.$ Let us find the form $\phi$(\textit{z}) using this initial condition. We replace
\[z=e^{bt}{\left(\frac{e^{ju}-1}{e^{ju}-\frac{\mu}{b}}\right)}^{\frac{b}{b-\mu}},\] 
and get
\[e^{ju}=\frac{1-\frac{\mu}{b}{\left(e^{-bt}z\right)}^{\frac{b-\mu}{b}}}{1-{\left(e^{-bt}z\right)}^{\frac{b-\mu}{b}}}.\] 
Therefore 
\[H(u,0)=\phi (e^{ju})=e^{juN}=\phi(z)={\left(\frac{1-\frac{\mu}{b}z^{\frac{b-\mu}{b}}}{1-z^{\frac{b-\mu}{b}}}\right)}^N,\] 
and by \eqref{GrindEQ__4_} we can write the characteristic function \textit{H}(\textit{u},\textit{t}) as
\begin{equation} \label{GrindEQ__5_} 
H(u,t)={\left\{\frac{e^{ju}-\frac{\mu}{b}-\frac{\mu}{b}(e^{ju}-1){exp \{\ }(b-\mu)t\}}{e^{ju}-\frac{\mu}{b}-(e^{ju}-1){exp \{\ }(b-\mu)t\}}\right\}}^N. 
\end{equation} 
We put b = $\mu$ when substituted into the characteristic function \eqref{GrindEQ__5_} an uncertainty of the type 0/0 arises. To reveal such uncertainty let us use L'Hospital's rule. Using L'Hospital's rule we find the derivatives with respect to $\mu$ of the numerator and denominator
\[\mathop{lim}_{\mu \rightarrow b}\frac{f\left(\mu\right)}{g\left(\mu\right)}=\mathop{lim}_{\mu \rightarrow b}\frac{e^{ju}-\frac{\mu}{b}-\frac{\mu}{b}(e^{ju}-1){exp \{\ }(b-\mu)t\}}{e^{ju}-\frac{\mu}{b}-(e^{ju}-1){exp \{\ }(b-\mu)t\}}=\mathop{lim}_{\mu \rightarrow b}
\frac{f^\prime \left(\mu \right)}{g^\prime \left(\mu \right)}=\] 
\[=\mathop{lim}_{\mu \rightarrow b}\frac{1+(e^{ju}-1){exp \{\ }(b-\mu)t\}-\mu t(e^{ju}-1){exp \{\ }(b-\mu)t\}}{1-bt(e^{ju}-1)}=\frac{(e^{ju}-1)bt-e^{ju}}{(e^{ju}-1)bt-1}.\] 
For \textit{b}=$\mu$ the characteristic function \eqref{GrindEQ__5_} takes the form
\begin{equation} \label{GrindEQ__6_} 
H(u,t)={\left\{\frac{(e^{ju}-1)bt-e^{ju}}{(e^{ju}-1)bt-1}\right\}}^N. 
\end{equation} 
The characteristic function \textit{H}(\textit{u},\textit{t}) satisfies the equality
\[{\left.\frac{\partial H(u,t)}{\partial u}\right|}_{u=0}=jm(t),\] 
where \textit{m}(\textit{t}) is the average value of the number of occupied servers in the queuing system. From \eqref{GrindEQ__6_} it is easy to obtain
\[{\left.\frac{\partial H(u,t)}{\partial u}\right|}_{u=0}={\left.jM\{i(t)e^{jui}\}\right|}_{u=0}=jM\{i(t)\}=jm(t)=Nj{exp \{\ }(b-\mu)t\}.\] 
Then the average value \textit{m}(\textit{t}) of the number of occupied devices is
\begin{equation} \label{GrindEQ__7_} 
m(t)=N{exp \{\ }(b-\mu)t\}.                                                          
\end{equation} 

\section{Research of a system with an unlimited number of devices and with a Poisson incoming flow}

We consider a Markov queuing system with an unlimited number of devices whose input receives a non-stationary Poisson flow of intensity $\lambda$(\textit{t}). The form of the function $\lambda$(\textit{t})\textit{ }is defined below. The times for servicing various customers are stochastically independent and equally distributed. The service time distribution function is exponential with the parameter $\mu$.

The queuing system is called Markov since the process \textit{i}(\textit{t}) is a Markov chain with continuous time. Therefore, the considered system M(\textit{t}){\textbar}M{\textbar}$\mathrm{\infty}$ belongs to the class of Markov queuing models, and to study it one can apply the theory of Markov chains with continuous time and, first of all, compose a direct system of Kolmogorov differential equations for the probability distribution
\[P(i,t)=P\{i(t)=i\}.\] 

Obviously, the direct system of Kolmogorov differential equations has the form
\[\frac{\partial P(i,t)}{\partial t}=-(\lambda (t)+i \mu)P(i,t)+\lambda(t)P(i-1,t)+(i+1) \mu P(i+1,t).\] 
Let us define the characteristic function of the number of occupied servers
\[G(u,t)=Me^{jui(t)}=\sum^\infty_{i=0}{e^{jui}}P(i,t),\] 
where $j=\sqrt{-1}$ is the imaginary unit. Similarly to \eqref{GrindEQ__2_} we write down the equation for the characteristic function \textit{G}(\textit{u},\textit{t})
\[\frac{\partial G(u,t)}{\partial t}+\frac{\mu}{j}\left(1-e^{-ju}\right)\frac{\partial G(u,t)}{\partial u}=\lambda (t)(e^{ju}-1)G(u,t).\] 
The solution to this equation is similar to \eqref{GrindEQ__2_}. We write down the general solution of equation
\[G(u,t)=\psi \left(\left(e^{ju}-1\right)e^{-\mu t}\right){exp \left\{\left(e^{ju}-1\right)e^{-\mu t}\int^t_0{\mu (\tau)}e^{\mu \tau}d \tau\right\}\ },\] 
where $\psi$(\textit{.}) is an arbitrary differentiable function.

Suppose that at the initial instant of time \textit{t}=0 and there are \textit{N} devices in the system. Let us find the form $\psi$(\textit{.}) using this initial condition. We replace
\[y=\left(e^{ju}-1\right)e^{-\mu t},\] 
\[e^{ju}=ye^{\mu t}+1.\] 
Then 
\[G(u,0)=\psi (e^{ju})=e^{juN}=\psi (y)={\left(y+1\right)}^N,\] 
and we can rewrite the general solution for the characteristic function \textit{G}(\textit{u},\textit{t})
\begin{equation} \label{GrindEQ__8_} 
G\left(u,t\right)={\left\{1-e^{-\mu t}+e^{ju}e^{- \mu t}\right\}}^N{exp \left\{\left(e^{ju}-1\right)e^{-\mu t}\int^t_0{\lambda \left(\tau\right)e^{\mu \tau}d \tau}\right\}\ }.\ \ \ \  
\end{equation} 
Find from \eqref{GrindEQ__8_}
\[{\left.\frac{\partial G(u,t)}{\partial u}\right|}_{u=0}={\left.jM\{i(t)e^{jui}\}\right|}_{u=0}=jM\{i(t)\}=jm(t)=jNe^{-\mu t}+je^{-\mu t}\int^t_0{\lambda (\tau)e^{\mu \tau}}d \tau,\] 
By virtue of this equality, the average value \textit{m}(\textit{t}) of the number of occupied devices has the form
\begin{equation} \label{GrindEQ__9_} 
m(t)=Ne^{-\mu t}+e^{-\mu t}\int^t_0{\lambda (\tau)e^{\mu \tau}d \tau}.                                                       
\end{equation} 

We choose the form of the function $\lambda$(\textit{t}) from the condition of coincidence of the mathematical expectations of the number of occupied devices in the servicing systems under consideration. Then, by virtue of \eqref{GrindEQ__7_} and \eqref{GrindEQ__9_}, we can write equation
\[N{exp \{\ }(b-\mu)t\}=Ne^{-\mu t}+e^{- \mu t}\int^t_0{\lambda (\tau)e^{\mu \tau}d \tau},\] 
from which we obtain that the intensity $\lambda$(\textit{t}) has the form
\begin{equation} \label{GrindEQ__10_} 
\lambda(t)=bN{exp \{\ }(b-\mu)t\}.                                                          
\end{equation} 
For such an intensity of the incoming Poisson flow in the system M(\textit{t}){\textbar}M{\textbar}$\mathrm{\infty}$  and the autonomous system the average values of the number of occupied devices are the same and are determined by the equality \eqref{GrindEQ__7_}. For the intensity $\lambda$(\textit{t}) of the form \eqref{GrindEQ__10_}, the characteristic function \eqref{GrindEQ__8_} takes the form
\begin{equation} \label{GrindEQ__11_} 
G(u,t)={\left\{1-e^{-\mu t}+e^{ju}e^{-\mu t}\right\}}^N{exp \left\{\left(e^{ju}-1\right)Ne^{- \mu t}\left(e^{bt}-1\right)\right\}\ }.                     
\end{equation} 

Note that the characteristic function of the sum of independent random variables is equal to the product of their characteristic functions. For a Poisson random variable with parameter \textit{a}, the characteristic function \textit{g}(\textit{u}) has the form
\begin{equation} \label{GrindEQ__12_} 
g(u)=\sum^ \infty_{k=0}{e^{juk}\frac{a^k}{k!}{exp \{\ }-a}\}={exp \{\ }-a\}\sum^\infty_{k=0}{\frac{(e^{ju}a)^k}{k!}}={exp \{\ }(e^{ju}-1)a\}.\ \ \  
\end{equation} 
For a random variable \textit{X} distributed binomially with parameters \textit{n }and \textit{p}, we can write
\[X=X_1+X_2+\dots X_n,\] 
where $X_i=\left\{ \begin{array}{c}
0,\ \mathrm{if\ event\ did}{\mathrm{n}}^{\mathrm{'}}\mathrm{t\ occure\ in\ i}\mathrm{-}\mathrm{th\ experience}, \\ 
1,\mathrm{if\ event\ occured\ in\ i}\mathrm{-}\mathrm{th\ experience}. \end{array}
\right.$

Since the \textit{X${}_{i}$ }are independent random variables according to the convolution theorem, the characteristic function \textit{h${}_{x}$}(\textit{u}) for a nonnegative random variable with binomial distribution has the form
\begin{equation} \label{GrindEQ__13_} 
h_X(u)=\prod^n_{i=1}{h_{X_i}}(u)=\prod^n_{i=1}{\left(1-p+pe^{ju}\right)}=(1-p+pe^{ju})^n.\ \ \ \ \ \ \ \ \ \ \ \ \ \  
\end{equation} 
It follows from \eqref{GrindEQ__11_} and \eqref{GrindEQ__12_}-\eqref{GrindEQ__13_} that \textit{G}(\textit{u},\textit{t}) is the characteristic function of the sum of two independent random variables, the first of which has a binomial distribution with probability of success $e^{-\mu t}$ and \textit{N -} the number of experiments, and the second has a Poisson distribution with a parameter $Ne^{- \mu t}\left(e^{bt}-1\right)$.

\section{Approximation of the birth process by the Poisson flow}
Let us find the probability distributions \textit{P}(\textit{i},\textit{t}) and \textit{PA}(\textit{i},\textit{t}), determined by the characteristic functions \textit{G}(\textit{u},\textit{t}) and \textit{H}(\textit{u},\textit{t}). Defining the inverse Fourier transform we obtain equalities
\begin{equation} \label{GrindEQ__14_} 
P(i,t)=\frac{1}{2 \pi}\int^\pi_{-\pi}{e^{-jui}G(u,t)du},                                                        
\end{equation} 

\begin{equation} \label{GrindEQ__15_} 
PA(i,t)=\frac{1}{2 \pi}\int^ \pi_{- \pi}{e^{-jui}H(u,t)du}.                                                       
\end{equation} 

Substituting \eqref{GrindEQ__5_} and \eqref{GrindEQ__11_} in \eqref{GrindEQ__14_} and \eqref{GrindEQ__15_}, and performing numerical integration here, we find the probability distributions \textit{P}(\textit{i},\textit{t}) and \textit{PA}(\textit{i},\textit{t}) for the given values of the parameters $\mu$, \textit{b}, \textit{N} and time \textit{t}. We define the distance between the distributions \textit{P}(\textit{i},\textit{t}) and \textit{PA}(\textit{i},\textit{t}) by the equation
\[\rho(b,t)=\mathop{max}_{0<i<\infty}\left|\sum^i_{n=0}{(P(n,t)-PA(n,t))}\right|.\] 

For the values of the parameters $\mu$=1 and \textit{N}=15 we find the distance values $\rho(b,t)$ for various values of the parameters \textit{b} (total fertility rate) and \textit{t} which are given in table. Here, by virtue of the equality $\mu$ =1, the unit of time is the average value of a person's lifespan (70-80 years old).
\begin{center}

\noindent \textbf{The distance between the distributions \textit{P}(\textit{i},\textit{t}) and \textit{PA}(\textit{i},\textit{t}) }

\noindent \textbf{for different values of the parameters \textit{b} and \textit{t}}

\noindent \textbf{}

\begin{tabular}{|p{0.2in}|p{0.4in}|p{0.4in}|p{0.4in}|p{0.4in}|p{0.4in}|p{0.4in}|} \hline 
\textit{t\newline b} & 0.1 & 0.2 & 0.4 & 0.6 & 0.8 & 1.0 \\ \hline 
0.8 & 0.009 & 0.019 & 0.039 & 0.057 & 0.075 & 0.090 \\ \hline 
1.2 & 0.014 & 0.030 & 0.059 & 0.086 & 0.112 & 0.136 \\ \hline 
1.5 & 0.019 & 0.038 & 0.073 & 0.107 & 0.125 & 0.126 \\ \hline 
1.6 & 0.020 & 0.040 & 0.079 & 0.105 & 0.108 & 0.096 \\ \hline 
1.7 & 0.021 & 0.042 & 0.082 & 0.096 & 0.086 & 0.068 \\ \hline 
1.8 & 0.023 & 0.045 & 0.080 & 0.082 & 0.065 & 0.045 \\ \hline 
1.9 & 0.024 & 0.048 & 0.076 & 0.067 & 0.046 & 0.039 \\ \hline 
\end{tabular}
\end{center}

The results given in table show that the distance $\rho(b,t)$ increases with increasing values of the parameter \textit{b} and time \textit{t}. 

\section{Conclusion}
Assuming that the approximation of one probability distribution to another is possible, when the discrepancy between them does not exceed 0.03, we can conclude that the approximation of the fertility process by the Poisson flow in the short-term forecasting is permissible, when \textit{t}$<$0.2. However, when considering the medium-term forecasts for 20-50 years, and even more long-term forecasts for periods of more than 50 years, the approximation of the birth process by the Poisson flow is inexpedient, since it leads to significant discrepancies.

\end{document}